\documentclass[a4paper,reqno,12pt]{amsart}

\begin{document}

\title[Kirkman's hypothesis revisited]
{Kirkman's hypothesis revisited}
\author[A.~Panholzer and H.~Prodinger]
{Alois Panholzer$^{\dag}$ and Helmut Prodinger$^{\ddag}$}

\address{${}^{\dag}$Institut f\"ur Algebra und Computermathematik\\
Technische Universit{\"a}t Wien\\ Wiedner Hauptstra{\ss}e 8--10\\ A-
1040 Wien\\ Austria.}
\email{Alois.Panholzer@tuwien.ac.at}\vskip0pt
\address{${}^{\ddag}$The John Knopfmacher
Centre for Applicable Analysis and Number Theory\\ Mathematics
Department\\ University of the Witwatersrand\\ P.O. Wits\\ 2050
Johannesburg\\ South Africa.}
\email{helmut@gauss.cam.wits.ac.za\vskip0pt {\it WWW-address{\rm:
}}\texttt{http://www.wits.ac.za/helmut/index.htm}}

\date{\today}

\begin{abstract}
Watson proved Kirkman's hypothesis (partially solved by Cayley).
Using  Lagrange Inversion, we drastically shorten Watson's
computations and generalize his results at the same time. 
\end{abstract}

\maketitle

Kirkman's hypothesis \cite{Kirkman1857} is (in changed notation)
the formula

\begin{align*}
   & \sum_{m=0}^{M} \sum_{n=0}^{N} \frac{1}{m+1} \binom{m+n}{n} \binom{2m+n+2}{m+n+2} \; 
   \times \\
   & \qquad \qquad  \times\frac{1}{M-m+1} \binom{M-m+N-n}{N-n} 
\binom{2(M-m)+N-n+2}{M-m+N-n+2} \\
   & =
   \frac{2}{M+2} \binom{M+N+1}{N} \binom{2M+N+4}{M+N+4} \, .
\end{align*}

Kirkman could not prove it, but Cayley \cite{Cayley1857} proved the special
case $N=0$ in 1857. After more than hundred years, Watson \cite{Watson62}
proved Kirkman's hypothesis by establishing the following
power series expansions. Set

\begin{equation*}   
   \psi(z,w):=\frac{1-w-2z-\sqrt{(1-w)^2-4z}}{2z(z+w)}\,,
\end{equation*}
then

\begin{align*}
\psi(z,w)&=\sum_{m,n}\frac{1}{m+1} \binom{m+n}{n} 
\binom{2m+n+2}{m+n+2} z^mw^n\,,\\
\psi^2(z,w)&=\sum_{m,n}
\frac{2}{m+2} \binom{m+n+1}{n} \binom{2m+n+4}{m+n+4} 
z^mw^n\,.
\end{align*}

Of course, Kirkman's hypothesis follows from this by writing
$\psi\cdot\psi=\psi^2$ and comparing coefficients.

However, Watson's derivation of these two expansions
required quite a bit of computation, in particular he
treated both cases differently and separately. 

Here, we present an extremely  simple computation using
the {\it Lagrange inversion formula\/} that has the advantage
of not only treating both cases together but rather finding the power
series expansion for $\psi^p(z,w)$ for general $p$. We refer for
the Lagrange inversion formula to \cite{GoJa83, Wilf94}.

The quadratic equation satisfied by $\psi(z,w)$ is
\begin{equation*}
   z(z+w)\psi^{2}(z,w)+(2z+w-1)\psi(z,w)+1=0 \, .
\end{equation*}

Writing $\psi = y/z$ and rearranging leads to the following equation of
Lagrange type:
\begin{equation*}
   y = z \frac{(1+y)^2}{1-w(1+y)} \, .
\end{equation*}

With the Lagrange inversion formula we obtain:
\begin{align*}
   [z^{m}w^{n}] \psi^{p}(z,w) & = [z^{m+p}w^{n}] y^{p}(z,w) = 
   \frac{p}{m+p}[y^{m}w^{n}] \left(\frac{(1+y)^{2}}{1-w(1+y)}\right)^{m+p} \\
   & = \frac{p}{m+p} \binom{m+n+p-1}{n} [y^{m}] (1+y)^{2m+n+2p} \\
   & = \frac{p}{m+p} \binom{m+n+p-1}{n} \binom{2m+n+2p}{m+n+2p} \, .
\end{align*}

This leads with $\psi^r\cdot\psi^s=\psi^p$ to the convolution formula 
{\it (generalized Kirkman
hypothesis)\/}:
\begin{align*}
   & \sum_{m=0}^{M} \sum_{n=0}^{N} \frac{r}{m+r} \binom{m+n+r-1}{n} \binom{2m+n+2r}{m+n+2r} \; 
   \times \\
   & \qquad  \qquad \times
\frac{s}{M-m+s} \binom{M-m+N-n+s-1}{N-n} \binom{2(M-m)+N-n+2s}{M-m+N-n+2s} \\
   & =
   \frac{p}{M+p} \binom{M+N+p-1}{N} \binom{2M+N+2p}{M+N+2p} \, .
\end{align*}

For other results of Kirkman's, treated with the Lagrange inversion formula, 
see \cite[Ex. 6.33-c]{Stanley99}.

\bibliographystyle{amsplain}


\providecommand{\bysame}{\leavevmode\hbox to3em{\hrulefill}\thinspace}

\end{document}